\input amstex
\documentstyle{amsppt}
\magnification = 1200
\hcorrection{.25in}
\def\({\left(}
\def\){\right)}
\def\[{\left[}
\def\]{\right]}

\def\R{9}

\NoBlackBoxes
\NoRunningHeads
\TagsOnRight
%\refstyle{C}
\topmatter
\title
Complete interpolating sequences for Paley-Wiener spaces and
Muckenhoupt's $(A_p)$ condition
\endtitle
\thanks{Research supported by NATO linkage grant LG 930329. While at MSRI,
Berkeley, during the fall 1995, the second author was supported in part by 
NSF grant DMS-9022140 and by the Norwegian Research Council.}
\endthanks
\author
Yurii I. Lyubarskii and Kristian Seip
\endauthor

\address
\hskip-\parindent
Yurii I. Lyubarskii,
Institute for Low Temperature Physics and Engineering, 47 Lenin Prospect,
310164 Kharkov, Ukraine
\endaddress
\email
lyubarskii\@ilt.kharkov.ua
\endemail
\address
\hskip-\parindent
Kristian Seip,
Department of Mathematical Sciences, The Norwegian Institute of Technology,
N--7034 Trondheim, Norway
\endaddress
\email
seip\@imf.unit.no
\endemail
%\dedicatory
%\enddedicatory
%\date
%\enddate
%\translator
%\endtranslator
%\keywords
%\endkeywords
\subjclass
30E05, 42A50
\endsubjclass
\abstract
We describe the complete interpolating sequences for the
Paley-Wiener spaces $L^p_\pi$ ($1<p<\infty$) in terms of Muckenhoupt's
$(A_p)$ condition. For $p=2$, this description coincides with those given
by Pavlov (1979),  Nikol'skii (1980), and Minkin (1992) of the unconditional 
bases of complex exponentials in $L^2(-\pi,\pi)$.  While the techniques of 
these authors are linked to the Hilbert space geometry
of $L^2_\pi$, our method of proof is based on turning the problem
into one about boundedness of the Hilbert transform in certain weighted $L^p$
spaces of functions and sequences.
\endabstract
\endtopmatter
\def\R{\Bbb R}
\def\Z{\Bbb Z}
\def\C{\Bbb C}
\def\Cpl{\C^+}
\def\Cmi{\C^-}

\def\L{\Lambda}
\def\l{\lambda}
\def\ol{\overline}
\def\ha{\frac12}
\def\eps{\varepsilon}

\def\H{\Cal H}
\def\Ca{\Cal C}
\document

\head 1. Introduction \endhead

In this paper we study interpolation in the Paley-Wiener spaces
$L^p_\pi$ $(1<p<\infty)$, which consist of all entire functions of
exponential type $\pi$ whose restrictions to the real line are in $L^p$. The
Paley-Wiener spaces are
Banach spaces when endowed with the natural $L^p(\R)$-norms. We want to
describe those sequences $\Lambda=\{\lambda_k\}$, $\lambda_k=\xi_k+i\eta_k$,
in the complex plane $\C$ for which the interpolation problem
$$
f(\lambda_k)=a_k
\tag{1}
$$
has a unique solution $f\in L^p_\pi$
for every sequence $\{a_k\}$ satisfying
$$
\sum_k|a_k|^pe^{-p\pi|\eta_k|}(1+|\eta_k|)<\infty.
\tag{2}
$$
Such sequences $\Lambda$ are termed {\it complete interpolating sequences
for} $L^p_\pi$. A classical example of a complete interpolating sequence for
$L^p_\pi$ ($1<p<\infty$) is the sequence of integers $ {\Bbb Z}$.

In the case $p=2$ this problem is equivalent
to that of describing all unconditional bases in $L^2(-\pi,\pi)$
of the form $\{\exp(i\l_kt)\}$. We refer to \cite{3} for an account of
this problem, including a detailed survey of its history.
The unconditional basis problem was solved by Pavlov  \cite{8}
under the additional restriction $\sup |\Im \l_k| <\infty$
and  by Nikol'skii \cite{7}, assuming only
$\inf \Im \l_k >-\infty$. Finally, Minkin \cite{6} solved the problem without
any a priori assumption on $\L$.

The methods of  \cite{3,6,7,8} are of a geometric
nature and make crucial use of the Hilbert space structure of $L^2_\pi$.
In this paper, we shall give a simpler proof, which works equally well
for all $p$, $1<p<\infty$. Incidentally, our method of proof  shows that
for $p=\infty$ or $0<p\leq 1$ there are no complete interpolating sequences.
(See also \cite{1}, which ``explains'' this curious phenomenon.)
The core of our approach is a careful study of properties of
the Hilbert transform
in weighted spaces of functions and its discrete version in weighted spaces of
sequences. More precisely, we turn our problem into one about boundedness
of the discrete Hilbert transform in a weighted space, defined on a
subsequence of $\L$ located in a horisontal strip, where the weight is
expressed in terms of certain infinite products involving
{\it all} the points of $\L$.

As an application of our main theorem, we prove a counterpart of the
well-known Kadets 1/4 theorem.
%complete interpolating sequences
%as well as related problems on bases of the form $\{\exp(i\l_kt)\}$ in the
%spaces $L^q(-\pi,\pi)$, $1<q<\infty$.

\head
2. Preliminary observations and statement of the main result
\endhead

%We will make use of some classical results about interpolation in the
%Hardy space
%$H^p$ of the half-plane, for which we refer to \cite{4}. For
%We also refer the reader
%to \cite{5}
%for the basic properties functions from the Paley-Wiener  spaces.
Suppose that $\L$ is a complete interpolating sequence for $L^p_\pi$.
Then $f\in L^p_\pi$ implies that $\exp(i\pi z)f(z)$ belongs
to the Hardy space $H^p$ of $\Cpl_a:=\{z\in\C:\ \Im z>a\}$ for each
$\ a\in \R$. It follows that the sequence $\L\cap\Cpl_a$
is $H^p$-interpolating in $\Cpl_a$ (see \cite{4}, Chapter 9). Similarly,
$\L\cap\Cmi_a$ is $H^p$-interpolating in the half-plane
$\Cmi_a:=\{z\in \C; \Im z<a\}$. So the sequences $\L\cap\Cpl_a$ and
$\L\cap\Cmi_a$ satisfy the Carleson condition in the corresponding
half-planes, i.e.,
$$
\sup_{\Im \l_j>a}\prod_{\Im \l_k>a, k\neq j}\left|\frac{\l_j-\l_k}
{\l_j-\overline{\l_k}-i2a}\right|>0, \quad
\sup_{\Im \l_j<a}\prod_{\Im \l_k<a, k\neq j}\left|\frac{\l_j-\l_k}
{\l_j-\overline{\l_k}-i2a}\right|>0.
\tag 3
$$
We note that, by standard manipulations (which we omit) with the Carleson
condition, this is equivalent to the following condition:
$$\sup_j \sum_{k,k\neq j}\frac{(1+|\eta_j|)(1+|\eta_k|)}{|\l_j-\l_k|^2}<\infty.
\tag 4$$
In particular, for some $\eps > 0$, the disks
$$
K(\l_k):=\{z: |z-\l_k|<10\eps (1+|\eta_k|)\}
$$
are pairwise disjoint.  (We fix this value of $\eps$
until the end of the paper.) Moreover, the measure
$$
\mu^+_\L:=\sum_{\eta_k>0} \eta_k\delta_{\l_k}
$$
($\delta_\l$ is the unit point measure at $\l$) is
a Carleson measure, i.e.,
$$
\int_{\Cpl} |f|^s d\mu^+_\L \leq C \|f\|_{H^s}^s
$$
for each function $f$ in the Hardy space $H^s(\Cpl)$,
$s\geq 1 $. Similarly, $\L$ generates a Carleson measure in
the lower half-plane as well as in each of the half-planes $\C^\pm_a$.

If $\L$ is a complete interpolating sequence for $L^p_\pi$, then
$$
\|f\|_{L^p(\R)} \leq C
\(\sum_k |f(\l_k)|^p e^{-p\pi|\eta_k|}(1+|\eta_k|)\)^{1/p},\quad
f\in L^p_\pi.
\tag{5}
$$
Indeed, since  $\L$ is interpolating
for $L^p_\pi$,  the operator
$$
T:f\mapsto \{f(\l_k)e^{-\pi|\eta_k|}(1+|\eta_k|)^{1/p}\}
$$
is bounded from $L^p_\pi$ onto $l^p$.  By the uniqueness of the solution of
the interpolation problem, we have $\ker T = \{0\}$, and it suffices to apply
the Banach theorem on inverse operators.

Given $x\in\R,\;r>0$, let  $Q(x,r)$ be the square
with center at $x$, side length $2r$, and sides parallel to the
coordinate axes. We say that a sequence $\L\subset \C$ is
{\it relatively dense} if there exists $r_0>0$ such that
$\L\cap Q(x,r_0)\neq \emptyset$ for each $x\in \R$.
If $\L$ is a complete interpolating sequence for $L^p_\pi$,
\thetag{5} forces $\L$ to be relatively dense: if this is not the case
and there exist  sequences $\{x_j\}\subset\R$ and $r_j\to \infty $ such that
$Q(x_j,r_j)\cap\L=\emptyset$, then, setting
$$
f_j(z)=\frac{\sin\frac{\pi}{2}(z-x_j)}{z-x_j},
$$
we get
$$
\sum_k |f_j(\l_k)|^pe^{-p\pi|\eta_k|}(1+|\eta_k|) \;\;\to\;\;0,\quad
                   j\to \infty,
$$
while $\|f_j\|_{L^p}$ is independent of $j$.

Suppose that $\L$ is a complete interpolating sequence for $L^p_\pi$.
Take $r>r_0$, where $r_0$ is as above, define
$$
Q_j=Q(4rj,r),\ \ \ \ j\in \Z,
$$
and pick a sequence $\Gamma=\{\gamma_j\}\subset \Lambda$ such that
$\gamma_j\in Q_j$. Let $\Sigma=\{\sigma_j\}$ be another sequence with
$|\gamma_j-\sigma_j|=\eps$. Suppose $w=\{w_j\}$ is a positive weight sequence.
Associate with it the weighted space $l^p_w$ consisting of
all sequences $a=\{a_k\}$ satisfying
$$
\|a\|_{w,p}^p:=\sum_k |a_k|^pw_k<\infty.
$$
We are interested in the boundedness of the {\it discrete Hilbert operator}
$\H_{\Gamma,\Sigma}$
defined by the relation
$$
\H_{\Gamma,\Sigma}:a=\{a_j\}\mapsto \{(\H_{\Gamma,\Sigma}a)_j\};\ \
(\H_{\Gamma,\Sigma}a)_j=\sum_k \frac{a_k}{\sigma_j-\gamma_k},
$$
on $l^p_w$. The following definitions are needed.
We say that $w$ satisfies the {\it discrete $(A_p)$ condition} if
$$
\sup_{k,n}\left(\frac{1}{n}\sum_{j=k+1}^{k+n} w_j\right)
\left(\frac{1}{n}\sum_{j=k+1}^{k+n} w_j^{-\frac{1}{p-1}}\right)^{p-1}<\infty.
$$
This condition is analogous to the classical continuous $(A_p)$ condition
for a positive weight $v(x)>0,\ x\in \R$:
$$
\sup_I \left\{\(\frac1{|I|}\int_I vdx\)
         \(\frac1{|I|}\int_Iv^{-1/(p-1)}dx\)^{p-1}\right\}<\infty ,
\tag{6}
$$
where $I$ ranges over all intervals in $\R$ (see \cite{2}). Recall that
the latter condition is necessary and sufficient for boundedness of the
classical Hilbert operator
$$
\H:f\mapsto (\H f)(t)=\frac1{i\pi}\int
              \frac{f(\tau)}{t-\tau}d\tau
$$
on the weighted space of functions $L^p(\R;v)$ consisting of all functions
$f$ satisfying
$$\|f\|_{v,p}^p:=\int |f(t)|^pv(t)dt <\infty.
$$
We shall need the following lemma.

\proclaim{Lemma 1} If $\H_{\Gamma,\Sigma}$ is bounded from
$l^p_w$ to $l^p_w$, then $w$
satisfies the discrete $(A_p)$ condition.
\endproclaim

\demo{Proof} We adopt the proof for the continuous case (see \cite{2}.)
Let $k$ and $n$ be given. For
convenience, put $I_1=\{k+1,k+2,...,k+n\}$, 
 $I_2=\{k+2n+1,k+2n+2,...,k+3n\}$.
Suppose that a positive sequence $a$ is supported on $I_1$. Then, for
$j\in I_2$, we have
$$
|(\H_{\Gamma,\Sigma}a)_j|
\geq\sum_l a_l\frac{\Re(\sigma_j-\gamma_l)}{|\sigma_j-\gamma_l|^2}
\geq \frac{C}{n}\sum_l a_l,
\tag{7}
$$
where $C$ is independent of $k$ and $n$.
Putting $a_l=1$, we get thus
$$
\sum_{j\in I_2} w_j \leq C \sum_{l\in I_1} w_l,
$$
and by symmetry
$$
\sum_{j\in I_1} w_j \asymp \sum_{l\in I_2} w_l;
\tag{8}
$$
here and in what follows the sign $\asymp$ means that the ratio of the
two sides lies between two positive constants.
Now we put $a_l=w_k^{\alpha}$ and get by \thetag{7}
$$
\left(\sum_{j\in I_2} w_j\right)
\left(\frac{1}{n}\sum_{l\in I_1} w_l^{\alpha}\right)^p\leq
C \sum_{m\in I_1}w_m^{1+\alpha p}.
$$
Finally, we put $\alpha=-\frac{1}{p-1}$ and invoke \thetag{8},
and the lemma is proved.
\qed
\enddemo

The converse of Lemma 1 is also true, but we will
not need that fact. Note also that the boundedness of the operator
$\H_{\Gamma,\Sigma}$ is independent of
the choice of sequence $\Sigma$, provided the condition
$|\gamma_j-\sigma_j|=\eps$ holds.

Let $\L$ be a complete interpolating sequence for $L^p_\pi$. If the
function $f_0\in L^p_\pi$ solves the interpolation problem
$f_0(\l_k)=\delta_{0,k},\ k\in\Z$, then $f_0(\mu)\neq 0$ for
$\mu \in \C\setminus\L$, since otherwise the function
$(z-\l_0)(z-\mu)^{-1}f_0(z)$ belongs to $L^p_\pi$ and vanishes on
$\L$. Since $f_0\in L^p_\pi$, $f_0$ belongs to the Cartwright class
$\Ca$ (see \cite{5}, Lecture 15) and, in particular,  the limit
$$
S(z)=\lim_{R\to\infty}\prod_{|\l_k|<R} (1-\frac{z}{\lambda_k})
\tag{9}
$$
exists and  defines the {\it generating function} of the sequence $\L$.
Besides, the solution $f_k\in L^p_\pi$ of the interpolation problem
$f_k(\l_n)=\delta_{k,n}$ has the form
$$
f_k(z)=\frac{S(z)}{S'(\l_k)(z-\l_k)}.
\tag{10}
$$

We may now formulate our main theorem.

\proclaim{Theorem 1} $\Lambda=\{\lambda_k\}$, where $\lambda_k=\xi_k+i\eta_k$,
is a complete interpolating sequence for $L^p_\pi$ if and only if the
following three conditions hold.
\roster
\item"{(i)}" The sequences $\Lambda\cap\C^+$ and $\Lambda\cap\C^-$
satisfy the Carleson condition in $\C^+$ and $\C^-$ respectively, i.e.
\thetag{3}
holds with $a=0$, and also $\inf_{k\neq j}|\l_k-\l_j|>0$.
\item"{(ii)}" The limit $S(z)$ in \thetag{9} exists and represents an
 entire function of exponential type $\pi$.
\item"{(iii)}" There exists a relatively dense subsequence
$\Gamma=\{\gamma_j\}\subset\L$ such that the sequence
$\{|S'(\gamma_j)|^p\}$ satisfies the discrete $(A_p)$ condition.
\endroster
Defining $F(x)=|S(x)|/\roman{dist}(x,\Lambda)$,  we may replace
statement $\roman{(iii)}$ by the following:
\roster
\item"{(iii')}"  $F^p$ satisfies the (continuous) $(A_p)$ condition.
\endroster
\endproclaim

Note that that condition (i) is equivalent to the statement that, for each
$a\in \R$, the sequences $\L\cap\C^\pm_a$ satisfy the Carleson condition (3).
Another, more compact way of expressing (i), is given by (4).

\head 3. Proof of Theorem 1: Necessity
\endhead

We have already proved the necessity of (i) and (ii), and also the existence
of a relatively dense sequence $\Gamma=\{\gamma_j\}\subset \L$.
We prove now that (iii) is necessary as well.
Let $\eps$ be as above. Then, for every $j$,
we can find a point $\sigma_j$ with $|\sigma_j-\gamma_j|=\eps$ and
$$
|S(\sigma_j)|=\eps |S'(\gamma_j)|.
$$
This follows from the fact that $S(z)(z-\gamma_j)^{-1}\neq 0$ for
$|z-\gamma_j|\leq\eps$, hence
$$
\min_{|z-\gamma_j|=\eps}|S(z)(z-\gamma_j)^{-1}|
   \leq |S'(\gamma_j)| \leq
     \max_{|z-\gamma_j|=\eps}|S(z)(z-\gamma_j)^{-1}|.
$$
Set $\Sigma=\{\sigma_j\}$. The
Plancherel-P\'{o}lya inequality (see \cite{5}, Lecture 20) yields
$$
\sum_j|f(\sigma_j)|^p\leq C \|f\|^p_{L^p},\ \ \ f\in L^p_\pi.
\tag{11}
$$
Now let $a=\{a_j\}$ be a finite sequence. By \thetag{10}, the unique solution
of the interpolation problem
$f(\gamma_j)=a_j$, $f(\lambda_k)=0$, $\lambda_j\not\in\Gamma$ has the
form
$$
f(z)=\sum_j \frac{a_j}{S'(\sigma_j)}\frac{S(z)}{(z-\gamma_j)}.
$$
By \thetag{5} and \thetag{11},  we have
$$
\sum_j |f(\sigma_j)|^p \leq C \sum_j |a_j|^p.
$$
Now, by our particular choice of the sequence  $\Sigma$, we obtain
(iii) by observing that Lemma 1 applies with $w_j=|S'(\gamma_j)|^p$.

To prove that (iii) implies (iii'), we need the following lemma.

\proclaim{Lemma 2}
Suppose $x\in \R$ and $\Re \gamma_j \leq x\leq \Re \gamma_{j+1}$.
Then there exists an $\alpha =\alpha(x)\in [0,1]$ such that
$$
|S'(\gamma_j)|^{\alpha} |S'(\gamma_{j+1})|^{1-\alpha}
\asymp |S(x)|/\roman{dist}(x,\Lambda),
$$
uniformly with respect to $x\in \R$. 
\endproclaim

In fact, assuming this lemma to hold, we see  that  (9) with 
$v=F^p$ follows from (iii) and the inequality
$t^\alpha s^{1-\alpha} \leq t+s, \ \ t,s>0,\ \ \alpha\in [0,1]$.

\demo{Proof of Lemma 2}
We assume that $x\in [\Re \gamma_j,\Re \gamma_{j+1}]$ and, for simplicity, 
$x\not\in\L$. Set $\L(x)=\{\l\in\L:|\l-x|<30r\}$.
(Here $r$ is the number used for  constructing  $\Gamma$.)
For $\alpha\in[0,1]$ we have
$$
\thickmuskip=.5 \thickmuskip
\medmuskip=.5 \medmuskip
\multline
\rho:= \frac{|S'(\gamma_j)|^\alpha |S'(\gamma_{j+1})|^{1-\alpha}}
{|S(x)|\text{dist}(x,\Lambda)^{-1}} \\
= \left\{\frac{\left |\frac1{\gamma_j}
                 \prod_{\l_k\in\L(x)\setminus\{\gamma_j\}}
		     \(1-\frac{\gamma_j}{\l_k}\)\right |^\alpha
              \left |\frac1{\gamma_{j+1}}
                 \prod_{\l_k\in\L(x)\setminus\{\gamma_{j+1}\}}
		     \(1-\frac{\gamma_{j+1}}{\l_k}\)\right |^{1-\alpha}}
{\left |\prod_{\l\in\L(x)} \(1-\frac x\l\)\right |}
      \text{dist}(x,\L)\right\} \\
\times\left\{\prod_{\l\in\L\setminus\L(x)}
       \frac{|\gamma_j-\l_k|^\alpha|\gamma_{j+1}-\l_k|^{1-\alpha}}
              {|x-\l_k|} \right \} =\Pi_1(x)\times\Pi_2(x).
\endmultline
$$
A simple estimation shows that $\Pi_1(x)\asymp 1$ uniformly with respect to 
$\alpha\in [0,1]$ so we need only estimate $\Pi_2(x)$.	      

Let us put
$$
\gamma_j=x-x_j+i y_j, \ \gamma_{j+1}=x+x_{j+1}+i y_{j+1}.
$$
The  values $x_j$ and $x_{j+1}$ depend upon $x$ and also satisfy
the inequalities $0\leq x_j, x_{j+1} \leq 8r$.
We may then write
$$
\multline
\rho^2\asymp \prod_{\l_k\not\in\L(x)}
\frac{((x-x_j-\xi_k)^2+(y_j-\eta_k)^2)^{\alpha} 
((x+x_{j+1}-\xi_k)^2+(y_{j+1}-\eta_k)^2)^{1-\alpha}}
{(x-\xi_k)^2+\eta_k^2} \\
=\prod_{\l_k\not\in\L(x)} \left(1-\frac{2x_j(x-\xi_k)+2y_j\eta_k+O(1)}
{(x-\xi_k)^2+\eta_k^2}\right)^{\alpha} \\
\times \left(1+\frac{2x_{j+1}(x-\xi_k)-2y_{j+1}\eta_k+O(1)}
{(x-\xi_k)^2+\eta_k^2}\right)^{1-\alpha}.
\endmultline
$$
Choosing $\alpha=\alpha(x)$ so that $\alpha x_j-(1-\alpha)x_{j+1}=0$, i.e.,
$\alpha=x_{j+1}/(x_j+x_{j+1})$,
we find that 
$$
\rho^2\asymp \exp\left(c\sum_{\l_k\not\in\L(x)}
                      \frac{|\eta_k|}{(x-\xi_k)^2+\eta_k^2}\right).
$$
By Carleson's condition \thetag{4}, 
the sum is uniformly bounded, and we are done.
\qed  
\enddemo

\head 4. Proof of Theorem 1: Sufficiency \endhead

We will now prove that (i), (ii), (iii') imply that $\Lambda$ is a complete
interpolating sequence.

To begin with, note that
$$
\int[ F(x)]^p\frac{dx}{1+|x|^p}<\infty
\tag{12}
$$
and
$$
\int [F(x)]^pdx=\infty.
\tag{13}
$$
The first relation follows from the fact that
$\int[F(x)]^p|\H f(x)|^pdx<\infty$ for each
bounded finite function $f$; it suffices to take
$f=\chi_{[0,1]}$. To obtain \thetag{12},
we may apply the operator $\H$ to a $\delta$-sequence $\{\delta_n(x)\}$.

First, we check that $\Lambda$ is a uniqueness set. To this end, we need
to estimate  $|S(z)|$ from below.

\proclaim{Lemma 3}
Let $\eps$ be the number from relation \thetag{3}. Then
$$
|S(z)|\geq C (1-|z|)^{-1/p} e^{\pi |\Im z|} \ \ \text{for} \ \ \roman{dist}
         (z,\L)>\eps(1+|\Im z|).
\tag{14}
$$
\endproclaim

\demo{Proof of Lemma 3}
Put $\L'=\L\cap\{z:|\Im z|<\eps\}$ and consider the auxiliary function
$$
S_1(z)=S(z)\prod_{\l\in\L'}\frac{z-\l+2i\eps}{z-\l}.
$$
It is plain that
$$
|S_1(z)|\asymp |S(z)|, \ \ |\Im z|> 3\eps,
\tag{15}
$$
and, besides, $|S_1(x)|^p$ satisfies the $(A_p)$ condition. Consider the
inner-outer factorization of $S_1$ in the upper half-plane,
$$
S_1(z)=e^{-i\pi z} G(z) B_1(z), \ \ \Im z>0.
\tag{16}
$$
Here the Blaschke product $B_1$ corresponds to the Carleson sequence
$(\L\cap\Cpl)\setminus\L'$ and, in particular,
$$
|B_1(z)|>c>0 \ \ \text{for dist}(z,\L)>\eps|\Im z|.
\tag{17}
$$
Moreover, $G$ is an outer function and $|G(x)|^p$ satisfies the $(A_p)$
condition. Therefore, $|G(x)|^{-q}$ is an $(A_q)$ weight (here  $1/p+1/q=1$),
$G(x)^{-1}(1+|x|)^{-1}\in L^q(\R)$, and thus
$$
\frac1{(z+i)G(z)}=\frac1{2\pi i} \int
    \frac1{(t+i)G(t)}\frac{dt}{t-z}, \ \ \Im z>0.
$$
It follows that
$$
\frac 1{|G(z)|} \leq C (1+|z|)^{1/p}.
\tag{18}
$$
Combining relations (15)--(18), we obtain (14) for $\Im z> 3\eps$.
The estimate for $\Im z<-3\eps$ is similar, and to fill the gap
$-3\eps <\Im z < 3\eps$, we may repeat the construction, taking
another horisontal line instead of $\R$.
\qed
\enddemo

Note that, since $|S_1(x)|^p$ is an $(A_p)$ weight, we have
$$
\int|S_1(x)|^pdx=\infty.
$$
The Phragm\'{e}n-Lindel\"{o}f theorem  (see \cite{5}, Lecture 20) yields
$$
\int|S_1(x+ia)|^pdx=\infty, \ \ a\in \R,
$$
and, by \thetag{15},
$$
\int |S(x+i)|^pdx =\infty.
$$
Again applying the Phragm\'{e}n-Lindel\"{o}f theorem, we get
$$
\int |S(x)|^pdx=\infty.
\tag{19}
$$

We are now in  position to prove the uniqueness. Indeed, if $f\in L^P_\pi$
and $f(\l)=0, \ \l\in\L$, then $\phi(z)=f(z)/S(z)$ is an
entire function of exponential type $0$,  and  \thetag{14}
yields that $|\phi(z)|$ is uniformly bounded for $z$ satisfying
$\text{dist}(z,\L)>\eps(|\Im z|+1)$. Therefore $\phi(z)\equiv C$,
which is incompatible with \thetag{19}, unless $C=0$.

It remains only to check that we can actually solve the interpolation
problem (1) for each sequence $a=\{a_k\}$ satisfying \thetag{2}.
It suffices to consider a finite sequence $a$ and bound
the norm of the solution by a constant times the left-hand side of \thetag{2}.
After doing so, we can apply a limit procedure.
If $a$  is a finite sequence, then, by (12),  the unique solution of
the interpolation problem (1) has the form
$$
f(z)=\sum_j \frac{a_k}{S'(\lambda_k)}\frac{S(z)}{(z-\lambda_k)}.
\tag{20}
$$
We split the sum (20) into two parts, corresponding to points lying
in $\Cpl\cup \R$ and in $\Cmi$, respectively. We may estimate the norm of
each sum
separately, so let us assume that all the $\lambda_k$ corresponding to
$a_k\neq 0$ are in
$\Cpl\cup\R$.
Clearly, we may estimate the $L^p$ integral along
$\Im(z)=-\ha$. Let us, however, for conventional reasons, estimate it along
$\R$ and assume all the points $\lambda_k$ satisfy $\eta_k\geq \ha$. Now
let
$$
B(z)=\prod_k\frac{z+\frac{i}{2}-(\lambda_k+\frac{i}{2})}
                 {z+\frac{i}{2}-(\ol{\lambda_k}-\frac{i}{2})}.
$$
Writing $S(z)=B(z)e^{-i\pi z} G(z)$, where $G$ is an outer function in $\Cpl$,
we observe that (iii') is equivalent to
$|G(x)|^p$ satisfying the $(A_p)$ condition. We note that we have
$$
|S'(\lambda_k)|\asymp |G(\lambda_k)|\frac{e^{\pi\eta_k}}{\eta_k}.
$$
Thus it is enough to consider the $L^p$ boundedness of
$$
\tilde{f}(x)=\sum_k\frac{a_k \eta_k e^{-\pi\eta_k}}{G(\lambda_k)}
\frac{G(x)}{x-\lambda_k}.
$$
By duality,
$$
\multline
\|\tilde{f}\|_p\asymp\sup_{\|h\|_q=1, h\in H^q}
                  \left|\sum_k\frac{a_k \eta_k e^{-\pi\eta_k}}
{G(\lambda_k)}\int_{\R} \frac{G(x)h(x)}{x-\lambda_k}dx\right| \\
\leq \sup_{\|h\|_q=1,h\in H^q}\left|\sum_k\frac{a_k \eta_k e^{-\pi\eta_k}}
{G(\lambda_k)}{\Cal H}Gh(\lambda_k)\right| \\
\leq \sup_{\|h\|_q=1,h\in H^q}
\left(\sum_k|a_k|^p \eta_k e^{-p\pi\eta_k}\right)^{1/p}
\left(\sum_k
\left|\frac{{\Cal H}Gh(\lambda_k)}{G(\lambda_k)}\right|^q\eta_k\right)^{1/q}.
\endmultline
$$
Since $|G(x)|^{-q}$ satisfies the $(A_q)$ condition, $G$ is an outer function
in $\Cpl$, and $h\in H^q, \|h\|_q\leq 1$, we have
${\Cal H}Gh(z)/G(z)\in H^q$, and  $\|{\Cal H}Gh(z)/G(z)\|_q\leq C$.
Since $\sum_k\eta_k\delta_{\lambda_k}$ is
a Carleson measure, we get the desired conclusion.

The sum corresponding to points in $\Cmi$ is treated similarly.

\head
5. A stability result
\endhead

We will now show how Theorem 1 can be used to obtain a result similar to
the Kadets 1/4 theorem. The same technique implies more sophisticated
stability results similar to the theorems of Avdonin and Katsnelson;
for these results we refer to \cite{3}.

For $1<p<\infty$ we denote by $q$ the conjugate exponent, $1/p+1/q=1$, and
put $$ p'=\max(p,q).$$ We may now prove:

\proclaim{Theorem 2}
Suppose that $\lambda_k=k+\delta_k$, $k\in {\Bbb Z}$. If
$|\delta_k|\leq d < 1/(2p')$
for every $k$, then $\L=\{\lambda_k\}$ is a complete interpolating sequence
for $L^p_\pi$. If merely $ |\delta_k|< 1/(2p')$
for every $k$, then $\L=\{\lambda_k\}$ is not necessarily a complete
interpolating sequence for $L^p_\pi$.
\endproclaim

Note that for $p=2$ this is precisely the Kadets theorem (see \cite{3}).

\demo{Proof of Theorem 2}
We prove first that the inequality $|\delta_k|<1/(2p')$ is not sufficient.
If $\delta_0=1$ and otherwise
$\delta_k=\text{sgn}(k)\delta$, $-1<\delta<1$, standard estimates of
infinite products yield
$$ F(x)\asymp (1+|x|)^{-2\delta}. $$
For $1<p<2$ we choose $\delta=1/(2q)$. Then
$$
\frac1{|x|}\int_0^{x} F^pdt
         \(\frac1{|x|}\int_0^x F^{-q}dt\)^{p-1}\geq \,C (\log(1+|x|))^{p-1},
$$
and the $(A_p)$ condition fails. We obtain the same conclusion if
$|\delta_k|<1/(2q)$ and
$\delta_k$ tends sufficiently fast to sgn$(k)/(2q)$ as $k$ tends
to $\pm\infty$. If $2<p<\infty$, we put $\delta=-1/(2p)$, and argue similarly.

With  $\Lambda$  as above, define $\lambda_\alpha=(k+\alpha \delta_k)$
and $\Lambda_\alpha=\{\lambda_\alpha\}$, where
$\alpha$ is a real number. Suppose that $\delta<1/2$ and $|\alpha|\delta<1/2$,
so that the distance between any two distinct points of $\Lambda$, and likewise
the distance between any two distinct numbers of $\L_\alpha$, exceeds a
certain positive number. Then estimates of infinite products show that
$$ F_\alpha(x)\asymp [F(x)]^\alpha, \tag{21}$$
where $F_\alpha(x)=|S_\alpha(x)|/\text{dist}(x,\Lambda_\alpha)$ and $S_\alpha$
is the generating function of $\L_\alpha$.

Suppose first that $1<p<2$. If
$d<1/(2q)$, then $F^2_{q/2}$ satisfies the $(A_2)$ condition, according to
the classical 1/4 theorem. By (21), it means that $F^q$ satisfies the
$(A_2)$ condition, which implies, by H\"{o}lder's inequality, that $F^p$
satisfies the $(A_p)$ condition.

If $2<p<\infty$, put $\alpha=p/2$ and argue
similarly. \qed
\enddemo

\goodbreak

\Refs

\ref\no 1
\by C. Eoff
\paper The discrete nature of Paley-Wiener spaces
\jour Proc. Amer. Math. Soc. \vol 123 \yr 1995 \pages 505--512
\endref

\ref\no 2
\by R. Hunt, B. Muckenhoupt, and R. Wheeden
\paper Weighted norm inequalities for the conjugate function and Hilbert
transform
\jour Trans. Amer. Math. Soc.
\vol 176 \yr 1973 \pages 227--251
\endref

\ref\no 3 \by S. V. Khrushchev, N. K. Nikol'skii, and B. S. Pavlov
\book Unconditional bases of exponentials and reproducing kernels
\bookinfo in ``Complex Analysis and Spectral Theory'',
 Lecture Notes in Math. Vol. 864 \publ Springer-Verlag
\publaddr Berlin/Heidelberg \yr 1981 \pages 214--335
\endref

\ref\no 4 \by P. Koosis
\book Introduction to $H_p$ Spaces
\publ Cambridge University Press
\yr 1980
\endref

\ref\no 5 \by B. Levin
\book Lectures on Entire Functions
\publ American Mathematical Society
\yr 1996
\endref

\ref\no 6 \by A. M. Minkin
\paper Reflection of exponents, and unconditional bases of exponentials
\jour St. Petersburg Math. J.  \vol 3 \yr 1992 \pages 1043--1068
\endref

\ref \no 7 \by N. K. Nikol'skii
\paper Bases of exponentials and the values of reproducing kernels
\jour Dokl. Akad. Nauk SSSR \vol 252 \yr 1980 \pages 1316--1320
\transl \jour English transl. in Sov. Math. Dokl. \vol 21 \yr 1980
\endref

\ref\no 8 \by B. S. Pavlov
\paper Basicity of an exponential system and Muckenhoupt's condition
\jour Dokl. Akad. Nauk SSSR \vol 247 \yr 1979 \pages 37--40
\transl \jour English transl. in Sov. Math. Dokl. \vol 20 \yr 1979 
%\pages 655--659
\endref

\endRefs
\enddocument